\documentclass[reqno]{amsart}

\usepackage[utf8]{inputenc}
\usepackage[T1]{fontenc}
\usepackage{amsthm}
\usepackage{amsmath}
\usepackage{amssymb}
\usepackage[inline]{enumitem}
\usepackage[dvips]{graphicx}
\usepackage{color}
\usepackage{comment}
\usepackage{hyperref}
\usepackage{url}
\usepackage{stackrel}
\usepackage[left=3.5cm, right=3.5cm, paperheight=11.5in]{geometry}
\usepackage{fancyhdr}
\usepackage{mathrsfs}
\usepackage{stmaryrd}
\usepackage{soul}
\usepackage{datetime}
\longdate

\newtheorem{theorem}{Theorem}
\newtheorem*{main*}{Main Theorem}
\newtheorem{lemma}{Lemma}
\newtheorem{conjecture}{\textbf{Conjecture}}

\theoremstyle{definition}

\newtheorem*{def*}{Definition}

\theoremstyle{remark}
\newtheorem{claim}{\textsc{Claim}}
\newtheorem*{claim*}{\textsc{Claim}}

\hypersetup{
    pdftitle={A characterization of Sophie Germain primes},
    pdfauthor={Paolo Leonetti},
    pdfmenubar=false,
    pdffitwindow=true,
    pdfstartview=FitH,
    colorlinks=true,
    linkcolor=blue,
    citecolor=green,
    urlcolor=cyan
}

\DeclareMathSymbol{\widehatsym}{\mathord}{largesymbols}{"62}

\renewcommand{\rho}{\varrho}

\begin{document}

\title{A characterization of Sophie Germain primes} 

\author{Paolo Leonetti}
\address{Universit\`a L. Bocconi, via Roentgen 1, 20136 Milano, Italy.}
\email{leonetti.paolo@gmail.com}



\begin{abstract}
\noindent{} Let $n\ge 5$ be an odd integer. It is shown that $\{1^{\sigma(1)},\ldots,n^{\sigma(n)}\}$ is a complete residue system modulo $n$ for some permutation $\sigma$ of $\{1,\ldots,n\}$ if and only if $\frac{1}{2}(n-1)$ is a Sophie Germain prime.   
Partial results are obtained also for the case $n$ even. 

\end{abstract}
\maketitle

\thispagestyle{empty}

\section{Introduction}\label{sec:introduction}

The aim of this article is to study an invariance property of complete residue systems modulo $n$, which turns out to be related to Sophie Germain primes. We recall that a prime $p$ is a \emph{Sophie Germain prime} if $2p+1$ is prime too, with the associated prime $2p+1$ which is then called a \emph{safe prime}. These special primes have applications in public key cryptography, pseudorandom number generation, and primality testing; see, for example, \cite{Agrawal, Matthews, Yap}. Originally, they have been used also in the investigation of cases of Fermat's last theorem \cite[\S{} 3.2]{Edwards}. It has been conjectured that there exist infinitely many Sophie Germain primes, but this remains unproven; cf., for instance, \cite[\S{} 5.5.5]{Shoup}.

Hereafter, we say that an integer $n\ge 2$ is \emph{nice} if $\{1^{\sigma(1)},\ldots,n^{\sigma(n)}\}$ is a complete residue system modulo $n$ for some permutation $\sigma$ of $\{1,\ldots,n\}$.  Then, our main result follows:
\begin{theorem}\label{th:main_odd}
Let $n\ge 5$ be an odd integer. Then $n$ is nice if and only if $n$ is a safe prime.
\end{theorem}

Partial results have been obtained also for the case $n$ even:
\begin{theorem}\label{th:main_even}
Let $n\ge 4$ be a nice even integer. Then $n=2p$ for some prime $p$ such that $p-1$ is squarefree. Conversely, if $n=2p$, for some safe prime $p\ge 7$, then $n$ is nice.
\end{theorem}


Note that, according to Theorem \ref{th:main_even}, $10$ is not a nice integer and, on the other hand, it is the double of a safe prime. However, the above results suggest the following:

\begin{conjecture}
An integer $n\ge 11$ is nice if and only if $n$ or $\frac{1}{2}n$ is a safe prime.
\end{conjecture}

Proofs of Theorem \ref{th:main_odd} and \ref{th:main_even} follow in \S\S{} \ref{sec:proofodd} and \ref{sec:proofeven}, respectively.

\subsection{Notations and conventions}
\label{sec:notations}
%
We let $\mathbf Z$ be the set of integers (endowed with its usual structure of ordered ring), $\mathbf N$ the non-negative integers, and $\mathbf N^+ = \mathbf N \setminus \{0\}$ the positive integers. Also, the set of (positive rational) primes $\{2,3,5,\ldots\}$ is denoted by $\mathbf{P}$.

Unless noted otherwise, the letters $n$, $m$, $i$, $j$, $k$, $t$ and $z$, with or without subscripts, will stand for positive integers, the letters $p$ and $q$ for primes, and the Greek letters $\sigma$ and $\eta$ for permutations.

Given an integer $n\ge 2$, we denote by $\mathbf{Z}_n$ the quotient ring between $\mathbf{Z}$ and its ideal $n\mathbf{Z}$; by an abuse of notation, sometimes we identify integers with its residue classes in $\mathbf{Z}_n$. The radical of $n$, that is, the product of the pairwise distinct primes which divide $n$, will be denoted by $\mathrm{rad}(n)$. Moreover, given $p \in \mathbf{P}$, the $p$-adic valuation of $n$ is $\upsilon_p(n)$, i.e., the greatest exponent $e \in \mathbf{N}$ for which $p^e$ divides $n$. 

Given integers $n,k\ge 2$, we denote by $\mathcal{A}_{n,k}$ the set of integers $m \in \{1,\ldots,n-1\}$ which are divisible by $k$, and by $\mathcal{Q}_{n,k}$ the set of (possibly zero) $k$-th power residues in $\mathbf{Z}_n$. The set of quadratic residues $\mathcal{Q}_{n,2}$ will be shortened with $\mathcal{Q}_n$.

Lastly, we write $\# S$ for the cardinality of a set $S$. We refer to \cite{Apo76} for basic aspects of number theory (including notation not defined here).

\section{Preliminaries}\label{sec:preliminaries}

Let us start settling down the cases of small values of $n$.
\begin{lemma}\label{lem:smallnice}
Every integer $n\in \{2,\ldots,7\}$ is nice.
\end{lemma}
\begin{proof}
It is enough to choose the permutation $\sigma$ according to the following table:
\begin{center}
\begin{tabular}{r|c|c|c|c|c|c|c|} 
\multicolumn{1}{r}{} &  \multicolumn{1}{c}{$\sigma(1)$}  & \multicolumn{1}{c}{$\sigma(2)$} & \multicolumn{1}{c}{$\sigma(3)$} & \multicolumn{1}{c}{$\sigma(4)$} & \multicolumn{1}{c}{$\sigma(5)$} & \multicolumn{1}{c}{$\sigma(6)$} & \multicolumn{1}{c}{$\sigma(7)$}\\ 
\cline{2-8}
$n=2$ & $1$ & $2$ & $\star$ & $\star$ & $\star$ & $\star$ & $\star$ \\ 
\cline{2-8} 
$n=3$ & $2$ & $1$ & $3$ & $\star$ & $\star$ & $\star$ & $\star$ \\ 
\cline{2-8}
$n=4$ & $2$ & $1$ & $3$ & $4$ & $\star$ & $\star$ & $\star$ \\ 
\cline{2-8}
$n=5$ & $2$ & $5$ & $1$ & $3$ & $4$ & $\star$ & $\star$ \\ 
\cline{2-8}
$n=6$ & $2$ & $1$ & $4$ & $5$ & $3$ & $6$ & $\star$ \\ 
\cline{2-8}
$n=7$ & $6$ & $2$ & $1$ & $5$ & $7$ & $3$ & $4$ \\ 
\cline{2-8}
\end{tabular}
\end{center}
\end{proof}

Accordingly, let us assume hereafter that $n\ge 8$.
\begin{lemma}\label{lem:firstsquarefree}
Let $n\ge 8$ be a nice integer. Then $n=r$ or $n=2r$ or $n=4r$ for some odd squarefree integer $r\ge 3$.
\end{lemma}
\begin{proof}
Let $\sigma$ be the associated permutation of $\{1,\ldots,n\}$. Note that $n$ divides $n^{\sigma(n)}$ and that $m^{\sigma(m)}$ is divisible by $\mathrm{rad}(n)$ for all $m \in \mathcal{A}_{n,\mathrm{rad}(n)}$. 

Since $n$ is nice by hypothesis, i.e., $\{1^{\sigma(1)},\ldots,n^{\sigma(n)}\}$ is a complete residue system in $\mathbf{Z}_n$, then $n$ does not divide $m^{\sigma(m)}$ for each $m \in \mathcal{A}_{n,\mathrm{rad}(n)}$. Moreover, since $\mathrm{rad}(n)$ divides $m$ by construction, then $n$ does not divide $\mathrm{rad}(n)^{\sigma(m)}$ for these integers $m$. In particular, $n$ does not divide $\mathrm{rad}(n)^{\#\mathcal{A}_{n,\mathrm{rad}(n)}}$. 

This implies that there exists $p \in \mathbf{P}$ which divides $n$ and 
$$
\upsilon_p(n)\ge 1+\#\mathcal{A}_{n,\mathrm{rad}(n)} = \frac{n}{\mathrm{rad}(n)} =\prod_{q \in \mathbf{P},\, q\mid n}q^{\upsilon_q(n)-1} \ge p^{\upsilon_p(n)-1}.
$$

If $p=2$, it follows that $\upsilon_2(n)=1$ or $\upsilon_2(n)=2$, and $\upsilon_q(n)=1$ for all other primes $q$ which divide $n$. Lastly, if $p\ge 3$, then $\upsilon_q(n)=1$ for all primes $q$ which divide $n$, i.e., $n$ is squarefree. 
\end{proof}

To conclude the section, we obtain a lower bound for the number of quadratic residues of a nice integer.
\begin{lemma}\label{lem:lowerboundQn}
Let $n\ge 2$ be a nice integer. Then $\# \mathcal{Q}_n \ge \lfloor \frac{1}{2}n\rfloor$.
\end{lemma}
\begin{proof}
Since $n$ is nice, the number of quadratic residues in $\{1^{\sigma(1)},\ldots,n^{\sigma(n)}\}$ has to be $\# \mathcal{Q}_n$. In particular, $\# \mathcal{Q}_n$ is greater than or equal to the number of even integers in $\{\sigma(1),\ldots,\sigma(n)\}$, that is, $\# \{1,\ldots,n\}\cap 2\mathbf{N}=\lfloor \frac{1}{2}n\rfloor$.
\end{proof}

\section{Proof of Theorem \ref{th:main_odd}}\label{sec:proofodd}

The proof will be splitted into two main parts.
\subsection{Only if part} Note that $5$ and $7$ are safe primes and, at the same time, are nice integers by Lemma \ref{lem:smallnice}. Hence, we can assume hereafter that $n$ is a nice odd integer $\ge 9$.

\begin{claim}\label{claim:nprime}
Let $n \ge 9$ be a nice odd integer. Then $n$ is prime.
\end{claim}
\begin{proof}
According to Lemma \ref{lem:firstsquarefree}, there are pairwise distinct odd primes $q_1,\ldots,q_k$ such that $n=q_1 \cdots q_k$. Note that, by the Chinese remainder theorem, the function $\mathbf{N}^+ \to \mathbf{N}^+$ defined by $n \mapsto \# \mathcal{Q}_n$ is multiplicative. Therefore, by Lemma \ref{lem:lowerboundQn}, we obtain
$$
\# \mathcal{Q}_1 \cdots \# \mathcal{Q}_k \ge \frac{1}{2}(q_1\cdots q_k-1),
$$
which simplifies to
$$
\prod_{i=1}^k \left(1+\frac{1}{q_i}\right) \ge 2^{k-1}\left(1-\frac{1}{q_1\cdots q_k}\right).
$$

Considering that $n\ge 8$, it follows that
$$
\left(1+\frac{1}{3}\right)\left(1+\frac{1}{5}\right)^{k-1}\ge \prod_{i=1}^k \left(1+\frac{1}{q_i}\right) \ge 2^{k-1} \left(1-\frac{1}{8}\right),
$$
which is satisfied only for $k=1$.
\end{proof}

Claim \ref{claim:nprime} will be refined further, by obtaining additional properties of nice primes.
\begin{claim}\label{claim:p1squarefree}
Let $p\ge 11$ be a nice prime. Then $p-1$ is squarefree.
\end{claim}
\begin{proof}
Let $\sigma$ be the associated permutation of $\{1,\ldots,p\}$ and suppose, for the sake of contradiction, that there exists a prime $q$ such that $q^2$ divides $p-1$. Then, note that if $m$ is a $q$-th power residue or if $\sigma(m)$ is divisible by $q$, then $m^{\sigma(m)}$ is a $q$-th power in $\mathbf{Z}_p$.

Since the number of $q$-th powers in $\{1^{\sigma(1)},\ldots,p^{\sigma(p)}\}$ has to be $\# \mathcal{Q}_{p,q}$ and $\# \mathcal{A}_{p,q}$ is smaller than $\# \mathcal{Q}_{p,q}$, it follows that $m \in \mathcal{Q}_{p,q}$ whenever $\sigma(m) \in \mathcal{A}_{p,q}$. In particular, $m^{\sigma(m)}$ is a $q^2$-th power in $\mathbf{Z}_p$. In turn, this implies that 
$$
\frac{p-1}{q}= \# \mathcal{A}_{p,q} \le \# \mathcal{Q}_{p,q^2}= 1+\frac{p-1}{q^2} \le 1+\frac{p-1}{2q}
$$

This is a contradiction because, on one hand, $q \ge \frac{1}{2}(p-1)$ by the above inequality, and, the other hand, $q \le \sqrt{p-1}$ by the fact that $q^2$ divides $p-1$.
\end{proof}

Without loss of generality, it can be assumed that, if $p$ is a nice (odd) prime with associated permutation $\sigma$, then 
\begin{equation}\label{eq:sigma1}
\sigma(1)=p-1.
\end{equation}
Indeed, by Fermat's little theorem, $m^{p-1}=1$ in $\mathbf{Z}_p$ for each $m \in \{1,\ldots,p-1\}$, implying that necessarily $\sigma(1)=p-1$ or $\sigma(p)=p-1$. On the other hand, if $p$ is a nice prime with associated permutation $\sigma$, then $p$ is a nice prime with another associated permutation $\tilde{\sigma}$ defined by $\tilde{\sigma}(p)=\sigma(1)$, $\tilde{\sigma}(1)=\sigma(p)$, and $\tilde{\sigma}(m)=\sigma(m)$ for each $m \in \{2,\ldots,p-1\}$.

To conclude the first part of the proof, it is enough to show the following:
\begin{claim}\label{claim:final1}
Let $p\ge 11$ be a nice prime. Then $p$ is a safe prime.
\end{claim}
\begin{proof}
Let $p\ge 11$ be a nice prime with associated permutation $\sigma$. According to Claim \ref{claim:p1squarefree}, $p-1$ is squarefree, i.e., there exist pairwise distinct odd $q_1,\ldots,q_k \in \mathbf{P}$ such that $p-1=2q_1\cdots q_k$ (note that $k\ge 1$ by the fact that $\frac{1}{2}(p-1) \ge 5$). Then, we claim that $k=1$.

Let us suppose, for the sake of contradiction, that $k\ge 2$ and define the (even) integers 
$$
z_1=\frac{p-1}{q_1}\text{ }\text{ and }\text{ }z_2=\frac{p-1}{q_{2}}.
$$
Then, at least one between $z_1$ and $z_2$ does not divide $\sigma(p)$. Indeed, in the opposite case, $p-1=\mathrm{lcm}(z_1,z_2)$ would divide $\sigma(p)$. On the other hand, since $\sigma(p)$ belongs to $\{1,\ldots,p\}$, then we have necessarily $\sigma(p)=p-1$, which contradicts \eqref{eq:sigma1}. Hence, there exists an integer in $\{z_1,z_2\}$, let us say $z$, which does not divide $\sigma(p)$, that is, $\sigma(p)$ does not belong to $\mathcal{A}_{p,z}$.

At this point, since $m^{\sigma(m)}$ is a $z$-th power in $\mathbf{Z}_p$ whenever $\sigma(m)$ belongs to $\mathcal{A}_{p,z}$ and $\# \mathcal{Q}_{p,z}=1+\#\mathcal{A}_{p,z}$, then $\mathcal{Q}_{p,z} \setminus \{0\}=\mathcal{A}_{p,z}$ in $\mathbf{Z}_p$. Denoting by $\xi$ a primitive root of $\mathbf{Z}_p$, it follows that there exists a permutation $\eta$ of $\{z,2z,\ldots,p-1\}$ such that 
$$
\{\xi^{z\eta(z)},\xi^{(2z)\eta(2z)},\ldots,\xi^{(p-1)\eta(p-1)}\}=\{\xi^z,\xi^{2z},\ldots,\xi^{p-1}\}
$$ 
in $\mathbf{Z}_p$. By Fermat's little theorem and the fact that $p$ is nice, we have by force $\eta(p-1)=p-1$. Therefore
$$
\{\xi^{z\eta(z)},\xi^{(2z)\eta(2z)},\ldots,\xi^{(p-1-z)\eta(p-1-z)}\}=\{\xi^z,\xi^{2z},\ldots,\xi^{p-1-z}\}
$$
in $\mathbf{Z}_p$, with the consequence that
$$
\{z\eta(z),(2z)\eta(2z),\ldots,(p-1-z)\eta(p-1-z)\}=\{z,2z,\ldots,p-1-z\}
$$
in $\mathbf{Z}_{p-1}$. Moreover, dividing all elements by $z$ and denoting by $q$ the prime $\frac{1}{z}(p-1)$, it follows that
$$
\{\eta(z),2\eta(2z),\ldots,(q-1)\eta(p-1-z)\}=\{1,2,\ldots,q-1\}
$$
in $\mathbf{Z}_q$. In particular, the products of the elements of each set must be the same in $\mathbf{Z}_q$. This is a contradiction, indeed the product of the set on the right is $(q-1)!\equiv -1\pmod{q}$ by Wilson's theorem, while on the left side
$$
(q-1)!\prod_{j=1}^{q-1}\eta(zj)=(q-1)!\prod_{j=1}^{q-1} zj=(q-1)!^2 z^{q-1} \equiv 1\pmod{q},
$$
by Fermat's little theorem and the fact that $\mathrm{gcd}(q,z)=1$.
\end{proof}


\subsection{If part}\label{subsectionif} Let $p\ge 5$ be a Sophie Germain prime. We claim that the prime $n=2p+1$ is nice.

Let $\xi$ and $\tau$ be generators of (the group of units of) $\mathbf{Z}_{n}$ and $\mathbf{Z}_{2p}$, respectively. Note that $\tau$ is odd. To conclude the proof of Theorem \ref{th:main_odd}, we have to construct an explicit permutation $\sigma$ of $\{1,\ldots,n\}$ such that $\{1^{\sigma(1)},\ldots,n^{\sigma(n)}\}$ is equal to $\{1,\ldots,n\}$ in $\mathbf{Z}_n$. 

To this aim, it is enough to set 
$$
\sigma(1)=2p, \text{ }\text{ }\text{ }\sigma(2p)=p,\text{ }\text{ }\text{ } \sigma(2p+1)=2p+1,
$$
together with
$$
\sigma\left(\xi^{(j \tau)^i}\right)=
\left\{
\begin{array}{l}
\!\! \left((j \tau)^i \ \ \ \ \bmod{2p}\right) \ \ \ \text{ if }\ i=0,\ldots,\frac{1}{2}(p-3) \\
\!\! \left((j \tau)^{i+1} \ \bmod{2p}\right)\ \ \ \text{ if }\ i=\frac{1}{2}(p-1),\ldots,p-2
\end{array}
\right.\!\!,
$$
for each $j \in \{1,2\}$, where $(x\bmod{2p})$ denotes the integer $y \in \{1,\ldots,2p\}$ such that $2p$ divides $x-y$.

Finally, let us check that this permutation really works. Define the sets
$$
\mathscr{A}_j=\{\xi^{(j\tau)^0},\xi^{(j\tau)^1},\ldots,\xi^{(j\tau)^{p-2}}\}
$$
in $\mathbf{Z}_{n}$, for each $j \in \{1,2\}$, and note that $\{1,2p,2p+1\} \cup \mathscr{A}_1 \cup \mathscr{A}_2$ is equal to $\{1,\ldots,n\}$ in $\mathbf{Z}_n$. Then, it is easy to see that, for each $j \in \{1,2\}$, the map 
$$
\mathscr{A}_j \to \mathscr{A}_j: m\mapsto m^{\sigma(m)}
$$
is actually a bijection. Indeed, for each $j \in \{1,2\}$, it holds
$$
\left(\xi^{(j \tau)^i}\right)^{\sigma\left(\xi^{(j \tau)^i}\right)}=
\left\{
\begin{array}{l}
\!\! \xi^{(j \tau)^{2i}} \ \ \ \ \ \ \ \text{ if }\ i=0,\ldots,\frac{1}{2}(p-3) \\
\!\! \xi^{(j \tau)^{2i+1}} \ \ \ \ \text{ if }\ i=\frac{1}{2}(p-1),\ldots,p-2 
\end{array}
\right.\!\!.
$$
This completes the proof. (Straightforward details are left to the reader.)


\section{Proof of Theorem \ref{th:main_even}}\label{sec:proofeven}

%
\subsection{First part} 

Note that $4$ and $6$ are nice integers by Lemma \ref{lem:smallnice} and both of them are in the form $2p$ for some prime $p$ such that $p-1$ is squarefree. Hence, let us hereafter that $n$ is a nice even integer $\ge 8$. In the same spirit of Claim \ref{claim:nprime}, we will prove that $\frac{1}{2}n \in \mathbf{P}$.

\begin{claim}\label{claim:upsiloneven}
Let $n \ge 8$ be a nice even integer. Then $n=2p$ for some prime $p$.
\end{claim}
\begin{proof}
According to Lemma \ref{lem:firstsquarefree}, there exist $\alpha \in \{1,2\}$ and pairwise distinct odd primes $q_1,\ldots,q_k$, with $k \ge 1$, such that $n=2^\alpha q_1 \cdots q_k$. Moreover, by Lemma \ref{lem:lowerboundQn} and the multiplicativity of $n\mapsto \#\mathcal{Q}_n$, we obtain
$$
\# \mathcal{Q}_{2^\alpha} \prod_{i=1}^k \# \mathcal{Q}_{q_i} \ge 2^{\alpha-1} q_1\cdots q_k.
$$

Considering that $\#\mathcal{Q}_{2^{\alpha}}=2$ for $\alpha \in \{1,2\}$ and $\#\mathcal{Q}_{q}=\frac{1}{2}(q+1)$ for each odd $q \in \mathbf{P}$, the above inequality simplifies to
$$
\prod_{i=1}^k \left(\frac{1}{2}+\frac{1}{2q_i}\right) \ge \frac{1}{2^{2-\alpha}}.
$$
On the other hand, note that, for all integers $k \ge 2$, it holds
$$
\prod_{i=1}^k \left(\frac{1}{2}+\frac{1}{2q_i}\right) \le \left(\frac{1}{2}+\frac{1}{2\cdot 3}\right)^k \le \left(\frac{2}{3}\right)^2 < \frac{1}{2}.
$$
It follows that $\alpha=1$ and $k=1$, i.e., $n=2p$ for some prime $p$.
\end{proof}

To complete the first part of the proof of Theorem \ref{th:main_even}, it will be enough to show that $p-1$ is squarefree. Accordingly, we will first show that $4$ does not divide $p-1$ and, then, that $q^2$ does not divide $p-1$ for each odd prime $q$.

Let $\sigma$ be a permutation associated to $2p$. Note that the number of quadratic residues in $\mathbf{Z}_{2p}$ is $p+1$, and, on the other hand, the number of even positive integers $\le 2p$ is $p$. It follows that $m$ has to be a quadratic residue whenever $\sigma(m)$ is even. Moreover, the residue modulo $2p $ of $m^{\sigma(m)}$ will be uniquely determined by the Chinese remainder theorem, given its residues modulo $p$ and modulo $2$ (in this respect, note that $m^k \equiv m\pmod{2}$ for all $m,k \in \mathbf{N}^+$).


\begin{claim}\label{claim:noresidue1mod4}
Let $p\ge 5$ be a prime such that $2p$ is nice. Then $4$ does not divide $p-1$.
\end{claim}
\begin{proof}
Let us suppose, for the sake of contradiction, that $4$ divides $p-1$. Then $p-1$ and $2p-1$ are quadratic residues in $\mathbf{Z}_{2p}$. By the above observations, at least one between $p-1$ and $2p-1$ has an even image under $\sigma$. This would contradict the fact that $1^{\sigma(1)}\equiv 1\pmod{2p}$ and $(p+1)^{\sigma(p+1)}\equiv p+1\pmod{2p}$ since, for all $k \in \mathbf{N}^+$, we have $(2p-1)^{2k}\equiv 1\pmod{2p}$ and $(p-1)^{2k} \equiv p+1 \pmod{2p}$.
\end{proof}

We conclude with the following:

\begin{claim}\label{claim:final}
Let $p\ge 7$ be a prime such that $2p$ is nice. Then $p-1$ is squarefree.
\end{claim}
\begin{proof}
Note $\frac{1}{2}(p-1)$ is odd by Claim \ref{claim:noresidue1mod4} and, by hypothesis, $\ge 3$. Hence, let us suppose, for the sake of contradiction, that there exists an odd prime $q$ such that $q^2$ divides $p-1$.

In addition, we have $\# \mathcal{Q}_{2p,q}=2+\frac{2}{q}(p-1)$ which is greater, on the other hand, than $\# \mathcal{A}_{2p,q}=\frac{2}{q}(p-1)$. With a reasoning similar to Claim \ref{claim:p1squarefree}, the number of $q^2$-th power residues in $\mathbf{Z}_{2p}$ has to be greater than or equal to the number of multiples of $q$ in $\{1,\ldots,2p\}$, implying that
$$
\frac{2(p-1)}{q}=\# \mathcal{A}_{2p,q} \le \# \mathcal{Q}_{2p,q^2}=2+\frac{p-1}{q^2} \le 2\left(1+\frac{p-1}{6q}\right).
$$

It follows that $q \ge \frac{5}{6}(p-1)$. This is a contradiction because $q\le \sqrt{p-1}$ by the fact that $q^2$ divides $p-1$ while, on the other hand, $\frac{5}{6}(p-1)>\sqrt{p-1}$ for all primes $p\ge 7$.
\end{proof}

\subsection{Second part}

Let $p\ge 7$ be a Sophie Germain prime. We claim that the integer $n=2(2p+1)$ is nice. The proof follows the same lines of reasoning in \S{} \ref{subsectionif}, therefore we provide here only a sketch.

Let $\xi$ and $\tau$ be generators of $\mathbf{Z}_{n}$ and $\mathbf{Z}_{2p}$, respectively. Then, define the permutation $\sigma$ of $\{1,\ldots,n\}$ by 
$$
\sigma(1)=2p, \sigma(2p)=p,\sigma(2p+1)=4p+1,\sigma(2p+2)=4p,\sigma(4p+1)=3p,\sigma(4p+2)=4p+2,
$$
together with
$$
\sigma\left((t\xi)^{(j \tau)^i}\right)=
\left\{
\begin{array}{l}
\!\! 2p(t-1)+\left((j \tau)^i \ \ \ \bmod{2p}\right) \ \ \ \text{ if }\ i=0,\ldots,\frac{1}{2}(p-3) \\
\!\! 2p(t-1)+\left((j \tau)^{i+1} \bmod{2p}\right)\ \ \ \text{ if }\ i=\frac{1}{2}(p-1),\ldots,p-2
\end{array}
\right.\!\!,
$$
for each $t,j \in \{1,2\}$, where $(x \bmod{2p})$ represents the integer $y \in \{1,\ldots,2p\}$ such that $2p$ divides $x-y$.

Finally, we have to check that this permutation really works. For each $t,j \in \{1,2\}$ define the sets
$$
\mathscr{A}_{t,j}=\{(t\xi)^{(j\tau)^0},(t\xi)^{(j\tau)^1},\ldots,(t\xi)^{(j\tau)^{p-2}}\}
$$
in $\mathbf{Z}_{n}$. Again, it is not difficult to check that, for each $t,j \in \{1,2\}$, the map 
$$
\mathscr{A}_{t,j} \to \mathscr{A}_{t,j}: m\mapsto m^{\sigma(m)}
$$
is actually a bijection. Indeed, for each $t,j \in \{1,2\}$, it holds
$$
\left((t\xi)^{(j \tau)^i}\right)^{\sigma\left(\xi^{(j \tau)^i}\right)}=
\left\{
\begin{array}{l}
\!\! (t\xi)^{(j \tau)^{2i}} \ \ \ \ \ \ \ \text{ if }\ i=0,\ldots,\frac{1}{2}(p-3) \\
\!\! (t\xi)^{(j \tau)^{2i+1}} \ \ \ \ \text{ if }\ i=\frac{1}{2}(p-1),\ldots,p-2 
\end{array}
\right.\!\!,
$$
which completes the proof. 


\section{Acknowledgements}
%
The author is grateful to Salvatore \textsc{Tringali} (University of Graz) for suggestions which improved the readability of the article.

\end{document}